\documentclass[12pt]{article}

\usepackage{amsfonts}
\usepackage{amsmath}
\usepackage{amssymb}
\usepackage{graphicx}
\usepackage[usenames]{color}

\usepackage{mathrsfs}

\newtheorem{theorem}{Theorem}[section]

\newtheorem{corollary}[theorem]{Corollary}

\newtheorem{remark}[theorem]{Remark}

\numberwithin{equation}{section} 

\def\en t{{{\rm Z}\mkern-5.5mu{\rm Z}}}

\def\<{\left<}
\def\>{\right>}
\def\({\left(}
\def\){\right)}

\def\9{{\infty}}
\def\barr{\begin{array}}
\def\earr{\end{array}}

\def\wt{\widetilde}
\def\wh{\widehat}

\def\lbb{{\lambda}}

\def\a{{\alpha}}

\def\n{\noindent }

\def\3{\subset }

\def\sk{\smallskip }

\def\bk{\bigskip }

\def\ve{{\varepsilon}}
\def\ora{\overrightarrow}

\def\bbe{{\mathbb{E}}}
\def\bbr{{\mathbb{R}}}
\def\bbp{{\mathbb{P}}}

\def\scrx{{\mathscr{X}}}
\def\scrp{{\mathscr{P}}}

\begin{document}

\begin{center}
{\Large{\bf Tridiagonal random matrix:

Gaussian fluctuations and deviations}}
\bigskip\bk

{\large{\bf Deng Zhang}}\footnote{Department of Mathematics,
Shanghai Jiao Tong University, 200240 Shanghai, China. This work is
supported by Notional Basic Research Program of China (973 Program,
2011CB808000, 2015CB856004) and Notional Natural Science Fundations
of China (11171216). } \footnote{Email address:
zhangdeng@amss.ac.cn}
\end{center}

\bk\bk\bk

\begin{quote}
\n{\small{\bf Abstract.} This paper is devoted to the Gaussian
fluctuations and deviations of the traces of tridiagonal random
matrix. Under quite general assumptions, we prove that the traces
are approximately normal distributed. Multi-dimensional central
limit theorem is also obtained here. These results have several
applications to various physical models and random matrix models,
such as the Anderson model, the random birth-death Markov kernel,
the random birth-death $Q$ matrix and the $\beta$-Hermite ensemble.
Furthermore, under some independent and identically distributed
condition, we also prove the large deviation principle as well as
the moderate deviation principle for the traces. } \\

\end{quote}

\vfill

\section{Introduction} \label{Intro}

We are here concerned with the tridiagonal random matrix
\begin{align} \label{TRM}
Q_n=&\left(
  \begin{array}{cccccc}
    d_1 & b_1 &   &   &   &   \\
    a_1 & d_2 & b_2 &  &   &   \\
      & a_2 & d_3 & b_3 &  &  \\
      &   & \ddots & \ddots & \ddots &   \\
    &   &   & a_{n-2}& d_{n-1} & b_{n-1} \\
      &   &  &   & a_{n-1} & d_n
  \end{array}
\right),
\end{align}
where $a_{i}, b_i$ and $d_i$, $1\leq i \leq n $, are random
variables with $a_0=b_n=0$.

Tridiagonal random matrix attracts significant interests in various
fields. In quantum mechanics, it is a finite-difference model of
one-dimensional random Schr\"{o}dinger operator, such as the
extensively studied Anderson model where $a_i=b_i =-1$ and all $d_i$
are independent and identically distributed (see e.g. \cite{PF92,
CL90}). The non-symmetric model with $a_i/b_i
>0$ also arises in the non-Hermitian quantum mechanics of Hatano and
Nelson, see e.g. \cite{GK05} and references therein. Tridiagonal
random matrix is also a basic model of random walks with random
environment in chain graphs, interesting examples include the random
birth-death Markov kernel proposed in \cite{BCC10, C09}, where
$a_{i-1} + d_i + b_i  =1$ and $\{(a_{i-1}, d_i, b_i)\}$ is an
ergodic random field, and the random birth-death $Q$ matrix recently
studied in \cite{HZ15}, where $d_i = -(a_{i-1} + b_i)$ and $\{a_i\}$
and $\{b_i\}$ are two sequences of strictly stationary ergodic
processes.

Moreover, tridiagonal random matrix also plays an important role in
random matrix theory. One well-known model is the $\beta$-Hermite
ensemble, in which $a_i=b_i$, $a_i$ is distributed as
$\beta^{-1/2}\chi_{i\beta}$ ($\chi_{i\beta}$ is the $\chi$
distribution with $i\beta$ degrees of freedom), $\beta>0$, $d_i$ is
normally distributed as $N(0,2/\beta)$, and $\{a_i,d_i\}$ are
independent. This model was first proposed in \cite{DE02} and
generalized the classical Gaussian ensembles for $\beta=1,2,4$,
corresponding to the Gaussian orthogonal ensemble (GOE), Gaussian
unitary ensemble (GUE) and Gaussian symplectic ensemble (GSE)
respectively. Due to its simple tridiagonal structure, one expects
to investigate some interesting spectral phenomena of general random
matrices from the study of the tridiagonal random matrix
\eqref{TRM}. For the generalization of the $\beta$-Hermite ensemble
to general symmetric tridiagonal random matrix (with independent
entries), we refer the interested reader to \cite{P09}, which is
another motivation for the present work.

The interest of this paper lies in the Gaussian fluctuations and
deviations of the trace of powers, or in other words, the moments of
empirical spectral distribution of the tridiagonal random matrix
\eqref{TRM}.

The fluctuations of the traces of random matrices are extensively
studied in the literature and, in a general context, turn out to be
Gaussian. For instance, the traces of the classical (unitary,
orthogonal and symplectic) compact groups were proved by Diaconis
and Shahshahani \cite{DS94}, by using the representation theory, to
be independent and normally distributed in the limit. For the Wigner
matrices, the Gaussian fluctuations of the traces are presented in
\cite[Theorem 2.1.31]{AGZ10} (see also \cite[Theorem 2.1]{G09}),
based on elaborate computations of moments. See also \cite{P09} for
the symmetric tridiagonal random matrix with indepenent entries,
where the Gaussian fluctuations of the traces were obtained by a
judicious counting of levels of paths. We also refer to
\cite{BS10,DE06,LP09,J98,So,So2} and references therein for
fluctuations of other linear eigenvalue statistics.

Furthermore, large deviation was well-known for the empirical
spectral distribution of the unitary invariant ensembles, including
the $\beta$-Hermite ensemble (see \cite{BeG97} and \cite{AGZ10}).
See also \cite{BeZ98} for the non self-ajoint matrix with
independent Gaussian entries. Moderate deviation principle for the
empirical spectral distribution of the Gaussian divisible matrix was
studied in \cite{DGZ03}.  Moreover, for the eigenvalue counting
function of the determinantal point process and Wigner matrix,  we
refer to \cite{DE13}
and \cite{DE13.2} . \\

Motivated by the physical models and random matrix models mentioned
above, we here consider the tridiagonal random matrix \eqref{TRM}
under quite general assumptions (see $(H.1)$ and $(H.2)$ in Section
\ref{CLT}, see also $(H.3)$ in Section \ref{LDP-MDP}), which, in
particular, allow the non-independent entries in \eqref{TRM}. We
prove the Gaussian fluctuations of the traces of such model.
Multi-dimensional central limit theorem is also given. Moreover, the
large deviations of the traces, in relation to those of additive
functionals of uniform Markov chains, as well as the moderate
deviation results are also obtained here. These results are
applicable to various models, such as the Anderson model, the random
birth-death Markov kernel, the random birth-death $Q$ matrix as well
as the symmetric tridiagonal random matrix, including the
$\beta$-Hermite ensemble.

Our proof is quite different from the standard method of moments
mentioned above and relies crucially on a path expansion of the
trace (see \eqref{trace-expan} below), which is formulated according
to the types of circuits determined by the tridiagonal structure of
\eqref{TRM}. This formula was recently used in our article
\cite{HZ15} to study the limiting spectral distribution of the
random birth-death $Q$ matrix. The advantage of this formula is
that, it reveals that the fluctuations and deviations of the traces
are asymptotically the same to those of the sum of $m$-dependent
random variables. This new point of view gives a unified way to
understand the limit behavior of traces of quite general tridiagonal
random matrices with non-symmetric structure or non-independent
entries. Moreover, it also allows to employ the analytic tools (e.g.
blocking arguments, Lyapunov's central limit theorem and large
deviations of additive functionals of Markov chains) to obtain the
Gaussian fluctuations and deviations of the traces. \\

The remainder of this paper is organized as follows. In Section
\ref{PRE} we first set up some preliminary notations and
definitions, then we present the path expansion formula
\eqref{trace-expan}. Section \ref{CLT} is devoted to the Gaussian
fluctuations of the traces, and Section \ref{LDP-MDP} is concerned
with the deviation results. Finally, Section \ref{DIS} includes some
discussions on the main results of this paper, and the Appendix,
i.e. Section \ref{APP} contains some technical proofs.\\

{\it Notations.} Throughout this paper, for $x\in\bbr^+$, $[x]$
denotes the largest integer not greater than $x$, $f =
\mathcal{O}(g)$ means that $|f/g|$ stays bounded, and $f_n=o(1)$
means that $|f_n|$ tends to zero, as $n\to \9$.

\section{Preliminaries} \label{PRE}

Let us first recall some notations from our recent paper
\cite{HZ15}. In particular, we will classify the circuits according
to their types. Then, we shall formulate the
path expansion formula of the trace precisely.\\

Let $k\geq 1$ be fixed. For every $0\leq l\leq [\frac k2]$, $1\leq
i\leq n-l$, set
\begin{align} \label{Qlmn}
   Q_{l,i}^{\overrightarrow{m}_l,\ora{n}_l}
   := \prod\limits_{j=0}^{l-1}
(a_{i+j}b_{i+j})^{m_{j+1}}  \prod \limits_{j=0}^l d_{i+j}^{n_j},
\end{align}
where, $\ora{m}_{l}=(m_1,m_2,...,m_l)$ and
$\ora{n}_{l}=(n_0,n_1,...,n_l)$. We say that $\ora{m}_l$ and
$\ora{n}_l$ are admissible, if $m_j\geq 0$ and $n_h\geq 0$ for every
$1\leq j\leq l$, $0\leq h\leq l$,
\begin{align} \label{Qlmn.2}
2\sum\limits_{j=1}^lm_j+\sum\limits_{h=0}^ln_h=k,
\end{align}
and if for some $l>0$, there exists $1\leq p\leq l$, such that
$m_p=0$, then $m_j=0$ and $n_h=0$ for all $ p<j\leq l$, $p\leq h\leq
l$. Set
\begin{align} \label{Psi}
\Psi_k := \{(l, \ora{m}_l, \ora{n}_l): 0\leq l\leq [\frac k2],
\ora{m}_l\ and\  \ora{n}_l\ are\ admissible\}.
\end{align}

The intuitions of these quantities can be seen as follows. By the
tridiagonal structure of \eqref{TRM}, we note that
\begin{align} \label{trace-expan.0}
  Tr Q_n^k
  =\sum\limits_{\pi\in \mathscr{C}_n} Q_{\pi},
  \ \ Q_{\pi}:=\prod\limits_{j=1}^k Q(\pi_j,\pi_{j+1}).
\end{align}
where $\mathscr{C}_n$ denotes the set of all circuits, i.e.
\begin{align*}
\mathscr{C}_n=\{\pi:(1,\cdots,k)\to
(1,\cdots,n):|\pi_j-\pi_{j+1}|\leq 1, 1\leq j\leq k,
\pi_{k+1}=\pi_1\}.
\end{align*}
The types of $\mathscr{C}_n$ can be determined by the set $\Psi_k$
and the vertices $i$, $1\leq i\leq n$. Indeed, $l$ is the largest
length between various vertices in the circuit $\pi$, and $i$ is the
leftmost vertex in $\pi$. The vector $\ora{m}_l$ determines the
circuit $\wh{\pi}$ consisting of the subindices of off-diagonal
entries in $Q_{\pi}$, and the coordinate $m_{j+1}$ is half of the
number of edges with the vertices $i+j$ and $i+j+1$. Similarly, the
vector $\ora{n}_l$ is related to the diagonal entries in $Q_{\pi}$,
and the coordinate $n_j$ is the number of loops with the vertex
$i+j$. Thus, $Q_{l,i}^{\overrightarrow{m}_l,\ora{n}_l}$ defined in
\eqref{Qlmn} can be viewed as a representative element of the type
$(l, \ora{m}_l,\ora{n}_l)$ and the leftmost point $i$.

Let $C_{l,i}^{\ora{m}_l,\ora{n}_l}$ denote the number of circuits of
the same type $(l,\ora{m}_l,\ora{n}_l)$ and the leftmost vertex $i$.
Note that for every $1\leq i,j\leq n$,
$C_{l,i}^{\ora{m}_l,\ora{n}_l}= C_{l,j}^{\ora{m}_l,\ora{n}_l}$. That
is, the number of circuits with the same type $\{l,\ora{m}_l,
\ora{n}_l\}$, though different leftmost vertices, are also the same.
Hence, we can set $C_{l}^{\ora{m}_l,\ora{n}_l}
:=C_{l,i}^{\ora{m}_l,\ora{n}_l}$.

With these notations, the expansion of the trace
\eqref{trace-expan.0} can be reformulated according to the types of
the circuits $\mathscr{C}_n$, i.e.
\begin{align} \label{trace-expan}
    Tr Q_n^{k}
   =\sum\limits_{(l,\ora{m}_l,\ora{n}_l)\in \Psi_k} C_l^{\ora{m}_l,\ora{n}_l}
\ \sum\limits_{i=1}^{n-l}Q_{l,i}^{\ora{m}_l,\ora{n}_l}.
\end{align}
This formula was recently used in our article \cite{HZ15} for the
study of limiting spectral distribution of the random birth-death
$Q$ matrix and, as we shall see later, is crucial for the
formulations and proofs of the main results in this paper.

We conclude this section by taking $Tr Q_n^3$ and $Tr Q_n^4$ for
examples. For  $Tr Q_n^3$, in this case $l=0,1$, $\Psi_k$ contains
the types $(0,(0),(3))$, $(1,(1),(1,0))$ and $(1,(1),(0,1))$, and
$C_l^{\ora{m}_l,\ora{n}_l}$ is equal to $1,3,3$ respectively. Hence,
by  \eqref{trace-expan} we have
\begin{align*}
  Tr Q_n^3
  =&\sum\limits_{i=1}^n d_i^3
     + 3\sum\limits_{i=1}^{n-1} (a_ib_i)d_i
     +3\sum\limits_{i=1}^{n-1} (a_ib_i)d_{i+1}.
\end{align*}
For $Tr Q_n^4$, $l=0,1,2$, $\Psi_k$ have the types $(0,(0),(4))$,
$(1,(1),(2,0))$, $(1,(1),(1,1))$, $(1,(1),(0,2))$, $(1,(2),(0,0))$
and $(1,(1,1),(0,0,0))$, and the corresponding
$C_l^{\ora{m}_l,\ora{n}_l}$ are $1,4,4,4,2$ and $4$ respectively.
Thus, it follows that
\begin{align*}
  Tr Q_n^4
  = \sum\limits_{i=1}^n d_i^4
    + \sum\limits_{i=1}^{n-1} a_ib_i(4d_i^2 + 4d_i d_{i+1} + 4d_{i+1}^2 + 2a_ib_i)
    +4\sum\limits_{i=1}^{n-2} a_ib_ia_{i+1}b_{i+1}.
\end{align*}

\section{Gaussian fluctuations}\label{CLT}

In this section, the main results Theorem \ref{THM-CLT-NONSTA} and
\ref{THM-CLT-NONSTA.2} are formulated in Subsection \ref{CLT-MAIN},
and then they are proved in Subsection \ref{PROOF-CLT} and
\ref{PROOF-CLT.2} respectively. Moreover, several applications are
also given in Subsection \ref{CLT-APP}.

\subsection{Main results} \label{CLT-MAIN}

Let $m$ be a fixed nonnegative integer. The random variables
$\{X_i\}$ are said to be $m$-dependent, if for any positive integers
$i$ and $j$ with $j-i>m$, $X_j$ is independent of the $\sigma$-field
generated by $\{X_h, 1\leq h\leq i\}$. In particular, $0$-dependence
is equivalent to independence.

Motivated by the physical models and random matrix models mentioned
in Section \ref{Intro}, we introduce the assumptions $(H.1)$ and
$(H.2)$ below.

\begin{enumerate}
\item[(H.1)] In the symmetric case (i.e. $a_i=b_i$ for all $i\geq 1$),
the off-diagonal entries $\{a_i\}$ are independent, and the diagonal
entries $\{d_i\}$ satisfy that $d_i= f(a_{i-1},a_i)$ with $f$ a
continuous function on $\bbr^2$ or \{$a_i,d_i$\} are all
independent.

In the non-symmetric case, the random vectors
$\{(a_{i-1},d_i,b_i)\}$ are independent.
\end{enumerate}

As regards the asymptotic behavior of the entries in \eqref{TRM}, we
assume that

\begin{enumerate}
\item[(H.2)] Let $\a$ and $\ve$ be two nonnegative constants, $0\leq \ve \leq
\a$. There exist constants $a, d, b$ and random variables $\eta_i,
\zeta_i, \xi_i$, $\eta, \zeta, \xi$, such that
\begin{align} \label{asym-abd}
   i^{-\a} (a_{i-1}, d_i, b_i)
   = (a,d,b) + i^{-\ve}(\eta_{i-1}, \zeta_i, \xi_i),
\end{align}
where $\eta_{i-1}, \zeta_i$ and $\xi_i$ satisfy
\begin{align} \label{asym-eatazeta}
   (\eta_{i-1}, \zeta_i, \xi_i) \overset{d}{\rightarrow}
   (\eta, \zeta, \xi),
\end{align}
$``\overset{d}{\rightarrow} "$ means convergence in distribution,
and for each $k\geq 1$,
\begin{align} \label{asym-eatazeta.2}
    \sup\limits_{i\geq 1} \bbe (\eta_i^{2k} + \zeta_i^{2k} + \xi_i^{2k}) <
    \9.
\end{align}
\end{enumerate}

\begin{remark} \label{Assumption-iid}
$(i)$. By H\"{o}lder's inequality, \eqref{asym-eatazeta.2} implies
that all moments of $\eta_n, \zeta_n, \xi_n$ are finite.

$(ii)$. Assumption $(H.1)$  allows to treat the tridiagonal random
matrix \eqref{TRM} with non-independent entries. The conditions on
$d_i$ in the symmetric case is mainly motivated by the random
birth-death Markov kernel and the random birth-death $Q$ matrix,
where $f(a_{i-1},a_i)$ is of the form $1-(a_{i-1}+a_i)$ and
$-(a_{i-1}+a_i)$ respectively.

$(iii)$. In the case that the entries of \eqref{TRM} are independent
and identically distributed (i.i.d.) and have all moments finite,
the assumptions $(H.1)$ and $(H.2)$ are obviously verified.

$(iv)$. In Assumption $(H.2)$ with $\a >0$, when $0<\ve<\a$ (resp.
$\ve =0$), the entries convergence in distribution to degenerate
(resp. non-degenerate) random variables. The degenerate case
actually includes the $\beta$-Hermite ensemble. See Subsection
\ref{CLT-APP} for more details.
\end{remark}

The main results in this section are formulated below. Taking into
account the complicated formulations of covariances in the case $\ve
>0$, we shall consider the case $\ve =0$ and $0<\ve \leq \a$ in
Theorem \ref{THM-CLT-NONSTA} and \ref{THM-CLT-NONSTA.2}
respectively, to make the structure of the variance and covariances
clean.

\begin{theorem}  \label{THM-CLT-NONSTA}
Assume $(H.1)$ and $(H.2)$ with $\ve =0$.

$(i)$. For each $k\geq 1$,  let $m_k=[\frac{k}{2}]$ (resp.
$[\frac{k}{2}]+1$) in the non-symmetric (resp. symmetric) case.  Set
$Tr \wt{Q_n^k} : = Tr Q_n^k - \bbe Tr Q_n^k$. Then,
\begin{align} \label{CLT-Nonsta}
   n^{-(\a k + \frac 12)} Tr \wt{Q_n^k}
   \overset{d}{\rightarrow} N(0, D_k).
\end{align}
Here
\begin{align} \label{A2-thm2}
   D_k= \frac{1}{2\a k +1} \bigg[Var (Z_{k,1}) + 2 \sum\limits_{j=1}^{m_k} Cov (Z_{k,1}, Z_{k,1+j}) \bigg],
\end{align}
and for each $1\leq i\leq m_k+1$,
\begin{align} \label{zi-thm2}
 Z_{k,i} = \sum\limits_{(l,\ora{m}_l,\ora{n}_l)\in \Psi_k} C_l^{\ora{m}_l,\ora{n}_l}
\prod\limits_{j=0}^{l-1} (\wt{\eta}_{i+j} \wt{\xi}_{i+j})^{m_{j+1}}
\prod\limits_{j=0}^{l} (\wt{\zeta}_{i+j})^{n_{j}},
\end{align}
where $\wt{\eta}_i, \wt{\zeta}_i, \wt{\xi}_i$ have the same
distributions as those of $a+\eta, d+\zeta$ and $b+\xi$
respectively, and they satisfy Assumption $(H.1)$ with all $a_i,
d_i, b_i$ replaced by $\wt{\eta}_i,\wt{\zeta}_i$ and $\wt{\xi}_i$
respectively.

$(ii)$ Given $k_1,\cdots, k_r \geq 1, r\geq 1$, set $m_{ij} = \max
\{m_{k_i}, m_{k_j}\}$. We have
\begin{align} \label{thm1-mult-clt}
  n^{-\frac 12}\ (n^{-\a k_1}Tr \wt{Q_n^{k_1}}, \cdots, n^{-\a k_r}Tr \wt{Q_n^{k_r}})
 \overset{d}{\longrightarrow } \Phi_{\Lambda},
\end{align}
where $\Phi_{\Lambda}$ is a $r$-dimensional normal distribution with
mean zero and covariance matrix $\Lambda$, for $1\leq i,j\leq r$,
\begin{align} \label{thm1-covar}
    \Lambda(i,j) =& \frac{1}{\a(k_i + k_j)+1}
     \bigg[Cov (Z_{k_i,1},Z_{k_j,1}) \nonumber \\
     & \qquad \qquad +
 \sum\limits_{h=1}^{m_{ij}} \(Cov(Z_{k_i,1},Z_{k_j,1+h}) +
Cov(Z_{k_i,1+h},Z_{k_j,1}) \) \bigg],
\end{align}
where $Z_{k_i,h}$ and $Z_{k_j,h}$, $1\leq h\leq m_{ij}+1$, are
defined as in \eqref{zi-thm2}.
\end{theorem}

\begin{theorem} \label{THM-CLT-NONSTA.2}
Assume $(H.1)$ and $(H.2)$ with $0< \ve \leq \a$. Let $Tr
\wt{Q_n^k}, m_k$ and $m_{ij}$ be as in Theorem \ref{THM-CLT-NONSTA}.

$(i)$. For each $k\geq 1$, $0\leq j \leq m_k$, there exists the
limit
\begin{align} \label{THM2-Var}
\sigma_k(1+j,1):= \lim\limits_{n\to \9} n^{-2(\a k - \ve)} Cov
(X_{k,n}, X_{k,n-j}),
\end{align}
and we have
\begin{align*}
    n^{-(\a k + \frac 12 - \ve)}
    Tr \wt{Q_n^k} \overset{d}{\to} N(0,D_k),
\end{align*}
where
\begin{align*}
   D_k = \frac {1}{2\a k   +1 -2 \ve} (\sigma_k(1,1)+
   2\sum\limits_{j=1}^{m_k}\sigma_k(1+j,1)).
\end{align*}

$(ii)$. Given $k_1,\cdots, k_r \geq 1$, $r\geq 1$. For every $1\leq
i,j\leq r$ and $0\leq h\leq m_{ij}$, there exists the limit
\begin{align} \label{THM2-Cov}
  \sigma_{k_i,k_j}(1+h,1):= \lim_{n\to \9} n^{-(\a (k_i + k_j) -2\ve)}
  Cov (X_{k_i,n}, X_{k_j, n-h}),
\end{align}
and we have
\begin{align*}
   n^{ - \frac 12 + \ve}
   (n^{-\a k_1} Tr \wt{Q_n^{k_1}}, \cdots, n^{-\a k_r}Tr \wt{Q_n^{k_r}} )
   \overset{d}{\rightarrow} \Phi_{\Lambda},
\end{align*}
where $\Phi_{\Lambda}$ is the $r$-dimensional normal distribution as
in Theorem \ref{THM-CLT-NONSTA}, but with the covariances
$\Lambda(i,j)$ defined by
\begin{align} \label{Cov-nonstat.2}
  \Lambda(i,j) =& \frac{1}{\a (k_i + k_j) +1 -2 \ve}
  \bigg[\sigma_{k_i,k_j}(1,1) \nonumber  \\
  &\qquad \qquad \qquad +\sum\limits_{h=1}^{m_{ij}}
   \(\sigma_{k_i,k_j}(1,1+h) + \sigma_{k_i,k_j}(1+h,1)\)\bigg],
\end{align}
where $\sigma_{k_i,k_j}(1,1+h)= \sigma_{k_j,k_i}(1+h,1)$.
\end{theorem}

\begin{remark}
Actually, in the case $0<\ve \leq \a$,  $\sigma_k(1+j,1)$ and
$\sigma_{k_i,k_j}(1+h,1)$ can be calculated explicitly in terms of
the covariances of $\eta_n, \xi_n$ and $\zeta_n$. Since the
formulations are complicated, we omit them in the statement of
Theorem \ref{THM-CLT-NONSTA.2}. Concrete calculations are shown in
Corollary \ref{COR-NONIID} below for the symmetric tridiagonal
random matrix motivated by \cite{P09} and the $\beta$-Hermite
ensemble.
\end{remark}

\subsection{Proof of Theorem  \ref{THM-CLT-NONSTA}} \label{PROOF-CLT}

The key observation for the proof is that, the path expansion
formula of the trace \eqref{trace-expan} indicates that the
fluctuation of $Tr Q_n^k$ is approximately the same as that of
\begin{align*}
     \sum\limits_{(l,\ora{m}_l,\ora{n}_l)\in \Psi_k} C_l^{\ora{m}_l,\ora{n}_l}
     \ \sum\limits_{i=1}^n  Q_{l,i}^{\ora{m}_l,\ora{n}_l}
   = \sum\limits_{i=1}^n X_{k,i},
\end{align*}
where
\begin{align} \label{Xi}
X_{k,i}
:=&\sum\limits_{(l,\ora{m}_l,\ora{n}_l)\in \Psi_k}
C_l^{\ora{m}_l,\ora{n}_l} Q_{l,i}^{\ora{m}_l,\ora{n}_l} \nonumber \\
=& \sum\limits_{(l,\ora{m}_l,\ora{n}_l)\in \Psi_k}
C_l^{\ora{m}_l,\ora{n}_l} \prod\limits_{j=0}^{l-1}
(a_{i+j}b_{i+j})^{m_{j+1}} \prod\limits_{j=0}^l d_{i+j}^{n_j}.
\end{align}
$\{X_{k,i}\}_{i\geq 1}$ is $m_k$-dependent with $m_k$ defined as in
Theorem \ref{THM-CLT-NONSTA}, due to Assumption $(H.1)$ and the
finite width of band in the tridiagonal random matrix \eqref{TRM}.
With this new point of view, we can expect the Gaussian fluctuation
of $Tr Q_n^k$, inspired by the fact that $m$-dependent random
variables are approximated normally distributed in the stationary
case (see e.g. \cite{HR48, C74, D05}). On the technical level, in
order to deal with the non-stationary case $\a>0$ in Assumption
$(H.2)$, we will employ the standard blocking arguments to separate
the sum $\sum_{i=1}^n
X_{k,i}$ into independent blocks and small ones. \\

{\it Proof of Theorem \ref{THM-CLT-NONSTA}.} $(i)$. We first note
that
\begin{align} \label{TrQ-SumX}
   \bbe \bigg| n^{-(\a k+\frac 12)} (Tr \wt{Q_n^k} - \sum\limits_{i=1}^n
   \wt{X}_{k,i}) \bigg|^2 \to 0,
\end{align}
where $\wt{X}_{k,i} = X_{k,i} - \bbe X_{k,i}$. (See the Appendix for
the proof.)

Thus, we only need to consider the fluctuation of $\sum_{i=1}^n
\wt{X}_{k,i}$, which is actually the sum of $m_k$-dependent random
variables.

Let $\kappa < \frac 14$, $n'=[n^{\kappa}]$, $p=[\frac {n}{n'}]$, $r=
n-pn'$. Set $\wt{Y}_{n,i}=n^{- \a k} \wt{X}_{k,i}$. For $1\leq i\leq
p$, let $\wt{U}_{n,i} := \sum\limits_{j=(i-1)n'+1}^{in'-m_k}
\wt{Y}_{n,j}.$ For $1\leq i\leq p-1$, set $\wt{Z}_{n,i} =
\sum\limits_{j=in'-m_k+1}^{in'} \wt{Y}_{n,j}$, $\wt{Z}_{n,p} =
\sum\limits_{j=pn'-m_k+1}^n \wt{Y}_{n,j}$ and $\wt{T}_n =
\sum\limits_{i=1}^p \wt{Z}_{n,i}$. Note that, $\{\wt{U}_{n,i}\}_{1
\leq i \leq p}$ are independent, and so are $\{\wt{Z}_{n,i}\}_{1\leq
i\leq p}$, due to the $m_k$-dependence of $\{X_{k,i}\}$. Moreover,
\begin{align*}
   n^{-(\a k + \frac 12 )} \sum\limits_{i=1}^n \wt{X}_{k,i}
   = n^{-\frac 12 } \sum\limits_{i=1}^n \wt{Y}_{n,i}
   = n^{-\frac 12 }\sum\limits_{i=1}^p \wt{U}_{n,i} + n^{-\frac 12 } \wt{T}_n.
\end{align*}

Let us show that, as $n\to \9$,
\begin{align} \label{tn-to-0}
   n^{- \frac 12 } \wt{T}_n \overset{d}{\longrightarrow} 0.
\end{align}
Indeed, by the independence of $\{\wt{Z}_{n,i}\}$,
\begin{align*}
    \bbe ( n^{- \frac 12 } \wt{T}_n)^2
   = \frac 1n \bbe (\sum\limits_{i=1}^p \wt{Z}_{n,i})^2
   = \frac 1n \sum\limits_{i=1}^p \bbe \wt{Z}_{n,i}^2.
\end{align*}
Then, it follows from the boundedness of moments in
\eqref{asym-eatazeta.2} that
\begin{align}  \label{sup-thm2-yni*}
   \bbe ( n^{- \frac 12 } \wt{T}_n)^2
   =  \mathcal{O}(\frac {p m_k^2}{n}) +  \mathcal{O}(\frac {(n-pn' +
m_k)^2}{n})
   \to 0,
\end{align}
which yields \eqref{tn-to-0}, as claimed. \\

In view of \eqref{tn-to-0} and  $pn'/n \to 1$, the proof of Theorem
\ref{THM-CLT-NONSTA} now reduces to proving that
\begin{align} \label{asym-thm2-uni}
   \frac {1}{\sqrt{pn'}} \sum\limits_{i=1}^p \wt{U}_{n,i}
\overset{d}{\rightarrow} N(0,D_k).
\end{align}

For this purpose, let us first treat the limit behavior of the
variance. Using the notations as in \cite{HR48}, we set
\begin{align} \label{def-ani}
A_{n,i}:= \bbe \wt{Y}_{n,i+m_k}^2 + 2 \sum\limits_{j=1}^{m_k} \bbe
\wt{Y}_{n,i+m_k} \wt{Y}_{n,i+m_k-j}.
\end{align}
By the independence of $\{\wt{U}_{n,i}\}_{1\leq i \leq p}$ and the
$m_k$-dependence of $\{\wt{Y}_{n,i}\}_{1\leq i\leq n}$,
straightforward computations show that
\begin{align} \label{exp-uni}
   & \bbe (\frac {1}{\sqrt{pn'}} \sum\limits_{i=1}^p \wt{U}_{n,i})^2
  = \frac {1}{pn'}  \sum\limits_{i=1}^p \bbe (\wt{U}_{n,i})^2 \nonumber  \\
  =& \frac {1}{pn'}\sum\limits_{i=1}^p
      \left[\bbe \(\sum\limits_{h=1}^{m_k} \wt{Y}_{n,(i-1)n'+h} \)^2
       + \sum\limits_{h=1}^{n'-2 m_k} A_{n,(i-1)n'+h} \right] \nonumber  \\
  =&  \mathcal{O}(\frac{m_k^2}{n'})
     + \frac {1}{pn'}\sum\limits_{i=1}^p
\sum\limits_{h=1}^{n'-2 m_k} A_{n,(i-1)n'+h},
\end{align}
where the last step is due to \eqref{asym-eatazeta.2}.

Set $n_{i,h}:= (i-1)n'+h+m_k$. In order to take the limit $p, n' \to
\9$, we note from  \eqref{def-ani} and the definition of
$\wt{Y}_{n,i}$ that
\begin{align} \label{cal-u}
    & \frac{1}{pn'} \sum\limits_{i=1}^p \sum\limits_{h=1}^{n'-2 m_k} A_{n,(i-1)n'+h} \nonumber \\
   =& \frac{1}{pn'} \sum\limits_{i=1}^p \sum\limits_{h=1}^{n'-2 m_k}
      \left[\bbe \wt{Y}^2_{n,n_{i,h}}
     + 2 \sum\limits_{j=1}^{m_k} \bbe \wt{Y}_{n,n_{i,h}} \wt{Y}_{n,n_{i,h}-j} \right]  \nonumber  \\
   =& \frac{1}{pn'} \sum\limits_{i=1}^p \sum\limits_{h=1}^{n'-2 m_k}
      (\frac{n_{i,h}}{n})^{2\a k} Var (\frac{X_{k,n_{i,h}}}{n_{i,h}^{\a
       k}})  \nonumber  \\
    &+ \frac{2}{pn'} \sum\limits_{j=1}^{m_k}
         \sum\limits_{i=1}^p \sum\limits_{h=1}^{n'-2 m_k}
        (\frac{n_{i,h}}{n})^{\a k} (\frac{n_{i,h}-j}{n})^{\a k}
    Cov (\frac{X_{k,n_{i,h}}}{n_{i,h}^{\a
       k}}, \frac{X_{k,n_{i,h}-j}}{(n_{i,h}-j)^{\a
       k}}).
\end{align}
It is not difficult to see that, for every $0\leq j\leq m_k$, as
$p,n' \to \9$,
\begin{align} \label{asym-ak}
      \frac {1}{pn'} \sum\limits_{i=1}^{p} \sum\limits_{h=1}^{n'-2 m_k}
      (\frac{n_{i,h}}{n})^{\a k} (\frac{n_{i,h}-j}{n})^{\a k}
      \to \frac{1}{2\a k+1}.
\end{align}
Moreover, we have that as $n(i,h)\to \9$,
\begin{align} \label{asym-cov}
     Cov (\frac{X_{k,n_{i,h}}}{n_{i,h}^{\a
       k}}, \frac{X_{k,n_{i,h}-j}}{(n_{i,h}-j)^{\a
       k}})
     \to Cov (Z_{k,1}, Z_{k,1+j}),
\end{align}
where $Z_{k,i}$, $1\leq i\leq m_k+1$, are defined as in
\eqref{zi-thm2}. (See the Appendix for the proof.)

Thus, plugging \eqref{asym-ak} and \eqref{asym-cov} into
\eqref{cal-u}, we obtain that as $p,n' \to \9$,
\begin{align} \label{asym-A}
  \frac{1}{pn'}\sum\limits_{i=1}^p
\sum\limits_{h=1}^{n'-2m_k} A_{n,(i-1)n'+h}
  \to  D_k
\end{align}
with $D_k$ defined as in \eqref{A2-thm2}, which together with
\eqref{exp-uni} yields that
\begin{align} \label{esti-2nd}
    \bbe (\frac{1}{\sqrt{pn'}}\sum\limits_{i=1}^p \wt{U}_{n,i})^2
   \to D_k.
\end{align}

As regards the third moments of $ \wt{U}_{n,i}/ \sqrt{n'}$ in
\eqref{asym-thm2-uni}, since $ \bbe |\wt{Y}_{n,j}|^3
=\mathcal{O}(1)$ and $n'=o(p^{\frac 13})$, we have for $1 \leq i
\leq p$,
\begin{align} \label{esti-3rd}
   \bbe \left|\frac{1}{\sqrt{n'}} \wt{U}_{n,i} \right|^3
   =& {(n')}^{-\frac 32} \bbe \left|\sum\limits_{j=(i-1)n'+1}^{in'- m_k}
\wt{Y}_{n,j} \right|^3 \nonumber \\
   =& \mathcal{O}({(n')}^{-\frac 32} (n'-m_k)^3)  \nonumber \\
    =&\mathcal{O}({(n')}^{\frac 32})
    =o(p^{\frac 12}).
\end{align}

Consequently, combining \eqref{esti-2nd} and \eqref{esti-3rd}, and
applying Lyapunov's central limit theorem (see
\cite[Appendix]{HR48}, see also \cite[Theorem 7.1.2]{C74}), we
obtain \eqref{asym-thm2-uni} and complete the proof of the first
part $(i)$.  \\

$(ii).$ In order to prove \eqref{thm1-mult-clt}, we just need to
show that for all $z_1,\cdots, z_r \in \bbr$, the linear combination
$n^{-1/2}\sum_{j=1}^r z_j n^{-\a k_j} Tr \wt{Q_n^{k_j}}$ is normally
distributed in the limit.

To this end, as in the proof of \eqref{TrQ-SumX}, the limit behavior
is the same as that of
\begin{align} \label{thm1-ii-linear}
    n^{-\frac 12}\sum\limits_{j=1}^r z_j n^{-\a k_j} \sum\limits_{i=1}^n
  X_{k_j,i}
     =n^{-\frac 12} \sum\limits_{i=1}^n \sum\limits_{j=1}^r z_j  n^{-\a k_j} X_{k_j,i},
\end{align}
where $X_{k_j,i}$ is defined as in \eqref{Xi} with $k_j$ replacing
$k$.

Since $\{X_{k_{j},i}\}_{i\geq 1}$ is $m_{k_j}$-dependent for each
$1\leq j\leq r$, $\{\sum_{j=1}^r z_j n^{-\a k_j} X_{k_j,i}\}_{i\geq
1}$ is $M$-dependent with $M=\max_{1\leq j\leq r} m_{k_j}$. Thus, as
in the preceding proof of the first assertion $(i)$, we deduce that
$n^{-1/2} \sum_{i=1}^n \sum_{j=1}^r z_j n^{-\a k_j} X_{k_j,i}$ is
normally distributed in the limit, which consequently implies the
Gaussian fluctuation of the random vectors $n^{-\frac 12} ( n^{-\a
k_1} Tr \wt{Q_n^{k_1}}, \cdots,  n^{-\a k_r} Tr
\wt{Q_n^{k_r}})$, $n\geq 1$. \\

It remains to compute the covariances $\Lambda(i,j)$. Let $m_{ij}=
\max \{m_{k_i}, m_{k_j} \}$. Then, $\{X_{k_i,q},X_{k_j,q}\}_{q\geq
1}$ is $m_{ij}$-dependent. Set $\wt{Y}_{k_i,n,q} := n^{-\a
k_i}\wt{X}_{k_i,q},$ and
\begin{align*}
   B_{n,q} :=& \bbe \wt{Y}_{k_i,n,q+m_{ij}} \wt{Y}_{k_j,n,q+m_{ij}} \\
   &+  \sum\limits_{h=1}^{m_{ij}} (\bbe \wt{Y}_{k_i,n,q+m_{ij}-h}  \wt{Y}_{k_j,n,q+m_{ij}}
      + \bbe \wt{Y}_{k_i,n,q+m_{ij}} \wt{Y}_{k_j,n,q+m_{ij}-h} ).
\end{align*}
Straightforward computations show that
\begin{align} \label{cal-cov-xij}
   &n^{-(\a (k_i +k_j) +1)} Cov (\sum\limits_{q=1}^n X_{k_i,q}, \sum\limits_{q=1}^n
   X_{k_j,q}) \nonumber \\
   =& n^{-1} \bbe (\sum\limits_{q=1}^{m_{ij}} \wt{Y}_{k_i,n,q})(\sum\limits_{q=1}^{m_{ij}} \wt{Y}_{k_j,n,q})
     + n^{-1}  \sum\limits_{q=1}^{n-m_{ij}} B_{n,q} \nonumber \\
   =& \mathcal{O}(n^{-1}) + n^{-1}  \sum\limits_{q=1}^{n-m_{ij}} \bbe \wt{Y}_{k_i,n,q+m_{ij}}
   \wt{Y}_{k_j,n,q+m_{ij}} \nonumber \\
    & + n^{-1} \sum\limits_{h=1}^{m_{ij}} \sum\limits_{q=1}^{n-m_{ij}}
      (\bbe\wt{Y}_{k_i,n,q+m_{ij}-h}  \wt{Y}_{k_j,n,q+m_{ij}}
      + \bbe\wt{Y}_{k_i,n,q+m_{ij}} \wt{Y}_{k_j,n,q+m_{ij}-h} ).
\end{align}
As in the proof of \eqref{asym-A}, we deduce that for each $0\leq h
\leq m_{ij}$, as $n \to \9$,
\begin{align*}
   \frac 1n \sum\limits_{q=1}^{n-m_{ij}}
      (\frac{q+m_{ij}-h}{n})^{\a k_i} (\frac{q+m_{ij}}{n})^{\a k_j}
   \to \frac{1}{\a (k_i+k_j) +1},
\end{align*}
and as $q \to \9$,
\begin{align*}
    Cov(\frac{X_{k_i,q+m_{ij}-h}}{(q+m_{ij}-h)^{\a k_i}}, \frac{X_{k_j,q+m_{ij}}}{(q+m_{ij})^{\a k_j}})
    \to Cov(Z_{k_i,1}, Z_{k_j,1+h}).
\end{align*}
Hence,
\begin{align} \label{asym-cov-kij}
   n^{-1} \sum\limits_{q=1}^{n-m_{ij}}
   \bbe (\wt{Y}_{k_i,n,q+m_{ij}-h} \wt{Y}_{k_j,n,q+m_{ij}})
   \to \frac{1}{\a (k_i+k_j) +1}Cov(Z_{k_i,1}, Z_{k_j,1+h}).
\end{align}
Similarly,
\begin{align} \label{asym-cov-kji}
   n^{-1} \sum\limits_{q=1}^{n-m_{ij}}
   \bbe (\wt{Y}_{k_i,n,q+m_{ij}}, \wt{Y}_{k_j,n,q+m_{ij}-h})
   \to \frac{1}{\a (k_i +k_j) +1}Cov(Z_{k_i,1+h}, Z_{k_j,1}).
\end{align}
Therefore, it follows from \eqref{cal-cov-xij}-\eqref{asym-cov-kji}
that, as $n\to \9$,
\begin{align*}
   n^{-(\a (k_i +k_j) +1)} Cov (\sum\limits_{q=1}^n X_{k_i,q}, \sum\limits_{q=1}^n X_{k_j,q})
   \to \Lambda(i,j).
\end{align*}
where $\Lambda(i,j)$ is defined as in \eqref{thm1-covar}, thereby
completing the proof. \hfill $\square$

\subsection{Proof of Theorem \ref{THM-CLT-NONSTA.2}} \label{PROOF-CLT.2}

{\it Proof of Theorem \ref{THM-CLT-NONSTA.2}.} The proof is similar
to that of Theorem \ref{THM-CLT-NONSTA}, we only need to prove that,
with the additional scaling $n^{2\ve}$,  the limits \eqref{THM2-Var}
and \eqref{THM2-Cov} exist. Below, we show the existence of the
limit \eqref{THM2-Var}, the proof of \eqref{THM2-Cov} follows
analogously.

Taking into account the definition of $X_{k,i}$ in \eqref{Xi}, we
need to  show the asymptotic estimate of $
\prod_{j=0}^{l-1}(a_{n+j}b_{n+j})^{m_{j+1}}
          \prod_{j=0}^{l} d_{n+j}^{n_j}.$

Indeed, by Assumption $(H.2)$, for each $0\leq j\leq l-1$,
\begin{align*}
    \frac{a_{n+j} b_{n+j}}{n^{2\a}}
    =& ab+ b\frac{\eta_{n+j}}{n^{\ve}} + a\frac{\xi_{n+j}}{n^{\ve}} +
\mathcal{O}(\frac{1}{n^{2\ve}}).
\end{align*}
Here, with a slight abuse of notation, $\mathcal{O}(n^{-2\ve})$
stands for the term of order $n^{-2\ve}$ after taking the
expectation. Hence, for $m_{j+1} \geq 1$,
\begin{align*}
    &(\frac{a_{n+j} b_{n+j}}{n^{2\a}})^{m_{j+1}} \\
   =& (ab)^{m_{j+1}} + (ab)^{m_{j+1}-1} m_{j+1} \(b\frac{\eta_{n+j}}{n^{\ve}} + a \frac{\xi_{n+j}}{n^{\ve}} \)
   +\mathcal{O}(\frac{1}{n^{2\ve}}).
\end{align*}
Hence, setting $|\ora{m}_l| = \sum\limits_{j=1}^l m_j$, we derive
that for $|\ora{m}_l| \geq 1$,
\begin{align} \label{calcu-ab-thm3}
    &\prod\limits_{j=0}^{l-1} (\frac{a_{n+j}
    b_{n+j}}{n^{2\a}})^{m_{j+1}} \nonumber \\
   =& (ab)^{|\ora{m}_l|-1} [ab +  a \sum\limits_{j=1}^l m_j
   \frac{\xi_{n+j-1}}{n^{\ve}}
    + b \sum\limits_{j=1}^l m_j \frac{\eta_{n+j-1}}{n^{\ve}} ]   +
\mathcal{O}(\frac{1}{n^{2\ve}}).
\end{align}

Similarly, setting $|\ora{n}_l| = \sum\limits_{j=0}^{l} n_j$, we
have that for $|\ora{n}_l| \geq 1$,
\begin{align} \label{calcu-d-thm3}
   \prod\limits_{j=0}^{l} (\frac{d_{n+j}}{n^{\a}})^{n_j}
   =& \prod\limits_{j=0}^{l}  \(d^{n_j} + d^{n_j-1} n_j \frac{\zeta_{n+j}}{n^{\ve}} \) +
   \mathcal{O}(\frac{1}{n^{2\ve}}) \nonumber \\
   =& d^{|\ora{n}_l|-1} [d +  \sum\limits_{j=0}^{l} n_j
   \frac{\zeta_{n+j}}{n^{\ve}} ]  + \mathcal{O}(\frac{1}{n^{2\ve}}).
\end{align}

Thus, \eqref{calcu-ab-thm3} and \eqref{calcu-d-thm3} yield that for
$|\ora{m}_l| \geq 1$, $|\ora{n}_l| \geq 1$,
\begin{align} \label{cal-abd.1-nonsym}
    &\prod\limits_{j=0}^{l-1} (\frac{a_{n+j} b_{n+j}}{n^{2\a}})^{m_{j+1}}
     \ \prod\limits_{j=0}^{l} (\frac{d_{n+j}}{n^{\a}})^{n_j} \nonumber  \\
   =& (ab)^{|\ora{m}_l|-1} d^{|\ora{n}_l|-1}
      \bigg[ abd + \sum\limits_{j=1}^l \bigg(adm_j \xi_{n+j-1}
    + bdm_j \eta_{n+j-1} \nonumber   \\
    &+ ab n_j
      \zeta_{n+j} \bigg) n^{-\ve}
      +ab n_0 \zeta_n n^{-\ve}  \bigg]
       + \mathcal{O}(n^{-2\ve}).
\end{align}

Moreover, it is easy to see that, for $|\ora{m}_l| = 0$,
\begin{align} \label{cal-abd.2-nonsym}
   (\frac{d_n}{n^{\a}})^k
   =d^{k-1} (d+ k \frac{\zeta_n}{n^{\ve}})
    + \mathcal{O}(\frac{1}{n^{2\ve}}),
\end{align}
and for $|\ora{n}_l| = 0$,
\begin{align} \label{cal-abd.3-nonsym}
    & \prod\limits_{j=0}^{l-1} (\frac{a_{n+j}
    b_{n+j}}{n^{2\a}})^{m_{j+1}} \nonumber \\
    =&(ab)^{\frac{k}{2}-1} [ab+ a\sum\limits_{j=1}^l m_j \frac{\xi_{n+j-1}}{n^{\ve}}
      + b\sum\limits_{j=1}^l m_j \frac{\eta_{n+j-1}}{n^{\ve}}]
      + \mathcal{O}(\frac{1}{n^{2\ve}}).
\end{align}

Plugging \eqref{Xi},
\eqref{cal-abd.1-nonsym}-\eqref{cal-abd.3-nonsym} into the right
hand side of \eqref{THM2-Var}, we consequently deduce that the limit
\eqref{THM2-Var} exists, which actually can be computed explicitly
in terms of the covariances of $\eta_n, \xi_n$ and $\zeta_n$. For
simplicity, we omit it here. The proof is complete. \hfill $\square$

\subsection{Applications} \label{CLT-APP}

Let us start with the independent identically distributed (i.i.d.)
case, which is a direct consequence of Theorem \ref{THM-CLT-NONSTA}
.
\begin{corollary} \label{COR-CLT-IID} ({\it The case $\a=0$}.)
In the symmetric case, assume that $\{a_i\}$ is a sequence of i.i.d.
random variables, $\{d_i\}$ satisfy that $d_i=f(a_{i-1},a_i)$ with
$f$ a continuous function, or $\{d_i\}$ are i.i.d random variables
and independent of $\{a_i\}$. In the non-symmetric case, assume that
$\{(a_{i-1},d_i,b_i)\}$ is a sequence of i.i.d random vectors.
Suppose also that in both cases the entries have all moments finite.
Then, the assertions in Theorem \ref{THM-CLT-NONSTA} hold. In
particular, the variance $D_k$ is given by
\begin{align} \label{Cor.1-Dk}
D_k= Var (X_{k,2}) + 2 \sum\limits_{j=1}^{m_k} Cov(X_{k,2},
X_{k,2+j}),
\end{align}
with $m_k$ defined as in Theorem \ref{THM-CLT-NONSTA} and $X_{k,i}$
as in \eqref{Xi}\footnote{The subindices of $X_{k,i}$ here starts
with $2$, since, e.g. in the case $d_i=f(a_{i-1},a_i)$, $i\geq 1$,
the distribution of $d_1(=f(0,a_1))$ is different from that of
$d_i=f(a_{i-1},a_i)$, $i\geq 2$. }, $2\leq i\leq m_k+2$, and the
covariances $\Lambda(i,j)$ are given by
\begin{align*}
\Lambda(i,j) = Cov (X_{k_i,2},X_{k_j,2}) +
 \sum\limits_{h=1}^{m_{ij}} (Cov(X_{k_i,2},X_{k_j,2+h}) +
Cov(X_{k_i,2+h},X_{k_j,2})),
\end{align*}
where $m_{ij}= \max\{m_{k_i},m_{k_j}\}$.
\end{corollary}

Corollary \ref{COR-CLT-IID} have several applications to the
physical models and random matrix models mentioned in Section
\ref{Intro}. \\

{\it Anderson model (\cite{PF92, CL90}).} When $a_i=b_i =-1$, and
$\{d_i\}$ are i.i.d random variables with all moments finite, this
tridiagonal random matrix \eqref{TRM} is referred to the Anderson
model in the physical literature, restricted to the bound domain
$\Lambda_n = \{1,\cdots, n\}$. Write $ Q_n = Q_{n,0} + V_n $, where
$V_n=diag(d_1,\cdots, d_n)$. $Q_{n,0}$ represents the finite
difference operator, a discrete analogue of the one dimensional
Schr\"{o}dinger operator, and $V_n$ represents the operator of
multiplication by the random field $d_i$, $1\leq i\leq n$.  For this
model, the conditions in Corollary \ref{COR-CLT-IID} are verified.
Hence, its traces of powers $Tr Q_n^k$, $k\geq 1$, are normally
distributed in the limit. \\

{\it Hatano-Nelson model (\cite{GK05}).} When $a_i/b_i >0$ and
$\{(a_{i-1},d_i,b_i)\}$ is a sequence of i.i.d. random vectors with
all moments finite, this tridiagonal random matrix is motivated by
the non-Hermitian quantum mechanics of Hatano and Nelson (see
\cite{GK05} and references therein). In this case, it follows from
Corollary \ref{COR-CLT-IID}
that the traces are approximated normally distributed. \\

{\it Random birth-death Markov kernel (\cite{BCC10, C09}).} This
tridiagonal random matrix arises in the random walks with random
environment in chain graphs. Given the chain graph $G=(V,E)$,
$V=\{1,\cdots,n\}$, $E=\{(i,j), |i-j|\leq 1 \}$. Consider the random
birth-death Markov kernel $Q_n$, with the state space $V$ and the
entries satisfying $a_i,b_i\in (0,1]$, $d_i\in [0,1]$, and
$a_{i-1}+d_i +b_i=1$. One concrete example is the random conductance
model, in which case $b_1=a_{n-1}=1$, $b_i=1-a_{i-1} =
U_{i,i+1}/(U_{i,i+1} + U_{i,i-1})$, and $\{U_{i,i+1}\}$ are i.i.d.
positive random variables. Another example is that $b_1=a_{n-1}=1$,
$b_i=1-a_{i-1} = V_i$, and $\{V_i\}$ are i.i.d. random variables on
$[0,1]$. We refer to \cite{BCC10} for more detail discussions and
for the study of the limiting spectral distribution of such model in
the ergodic environment. Here, in the i.i.d. environment, i.e.
$\{(a_{i-1}, d_i, b_i)\}$ is a sequence of i.i.d. random vectors
with all moments finite, we deduce from Corollary \ref{COR-CLT-IID}
that the traces are
asymptotically normally distributed. \\

{\it Random birth-death Q matrix (\cite{HZ15}).} This model is
motivated from the infinitesimal generator of continuous-time Markov
process in the chain graph. In this case, $a_i, b_i >0$ and $d_i=
-(a_{i-1}+ b_{i})$. In \cite{HZ15} we proved the existence of the
limiting spectral distribution in the strictly stationary ergodic
case (see \cite[Theorem 1]{HZ15}). Here, in the case that $\{a_i\}$
and $\{b_i\}$ are two sequences of i.i.d random variables, $\{a_i\}$
is independent of $\{b_i\}$ (or $a_i=b_i$, $i\geq 1$), and they have
all moments finite, Corollary \ref{COR-CLT-IID}
implies  the Gaussian fluctuations of the traces. \\

Now, we formulate and prove Corollary \ref{COR-NONIID} below for the
case $\a>0$, which is mainly motivated by the symmetric tridiagonal
random matrix studied in \cite{P09} and includes, in particular, the
well-known $\beta$-Hermite ensemble.

For simplicity, we will regard the $l+1$-dimensional vector
$\ora{n}_l$ in $\Psi_k$ in \eqref{Psi} as the vector in the larger
space $\bbr^{[\frac k2]+1}$, by just adding zeros to the remaining
$[\frac k2]-l$ coordinates, and for each $0\leq q\leq [\frac k2]$,
$\ora{e}_q$ denotes the $q$-th normal basis in
$\bbr^{[\frac{k}{2}]+1}$, namely, the $q$-th coordinate of
$\ora{e}_q$ is $1$ and the others are zeros.

\begin{corollary} \label{COR-NONIID} ({\it The case $\a>0$}.)
Consider the symmetric tridiagonal random matrix \eqref{TRM}, i.e.
$a_i=b_i$, $i\geq 1$.  Assume that in $(H.1)$ all $a_i$ and $d_i$
are independent. Assume also $(H.2)$ with $0< \ve \leq \a$ and
$\sup_{n\geq 1} \bbe |d_n|^k< \9$ for any $k\geq 1$. Then, the
assertions in Theorem \ref{THM-CLT-NONSTA.2} hold. In particular,
the covariances $\Lambda(i,j)$ are given by
\begin{align} \label{Cor-noniid-covariance}
   \Lambda(i,j) =
   \left\{
     \begin{array}{ll}
       \frac{a^{k_i+k_j-2}Var(\eta)}{\a(k_i+k_j)+1-2\ve} k_ik_j
       {k_i \choose k_i/2}{k_j \choose k_j/2}, & \hbox{if $k_i, k_j$ even;} \\
       \frac{a^{k_i+k_j-2}Var(\zeta)}{\a(k_i+k_j)+1-2\a} k_ik_j
       {k_i -1 \choose (k_i -1)/2}{k_j -1 \choose (k_j -1)/2}, & \hbox{if $\ve=\a$, $k_i, k_j$ odd ;} \\
       0, & \hbox{otherwise.}
     \end{array}
   \right.
\end{align}
\end{corollary}

\begin{remark} This result coincides with Theorem $3$ in
\cite{P09},\footnote{In \cite[(3.5)]{P09}, $m_2^{k+l}$,
$m_2^{\mathbb{I}_h(\gamma_1)+\mathbb{I}_h(\gamma_2)}$ and the
denominator $\a(k+l)$ shall be modified by $m_2^{\frac 12(k+l)}$,
$m_2^{\frac 12 (\mathbb{I}_i(\gamma_1)+\mathbb{I}_i(\gamma_2))}$ and
$\a(k+l)+1-2\a$ respectively.} under the conditions considered here.
However, the proof presented below is quite different from that in
\cite{P09}, it is analytic and mainly based on Theorem
\ref{THM-CLT-NONSTA.2}.
\end{remark}

{\it Proof.} First, by the condition on $d_n$, we note that in
Assumption $(H.2)$ $d=0$,  $\zeta_n= n^{-(\a-\ve)} d_n$, and the
limit $\zeta=0$ when $\ve < \a$. Below we shall discuss the case $k$
is even or odd respectively.

When $k$ is even, \eqref{Qlmn.2} indicates that
$|\ora{n}_l|:=\sum_{h=0}^l n_h$ is also even.  Since  $d=0$, from
\eqref{cal-abd.1-nonsym}-\eqref{cal-abd.3-nonsym} we see that the
main contribution comes from the case $|\ora{n}_l|=0$, i.e. there is
no loop in the circuit. Hence, it follows from \eqref{Xi} and
\eqref{cal-abd.3-nonsym} that
\begin{align} \label{Cor-X-even}
n^{-\a k} X_{k,n} = \sum\limits_{(l,\ora{m}_l,\ora{n}_l) \in \Psi_k}
 C_l^{\ora{m}_l, \ora{n}_l}
 [a^{k-2} (a^2 + 2a \sum\limits_{q=1}^l m_q
\frac{\eta_{n+q-1}}{n^{\ve}}) + \mathcal{O}(\frac{1}{n^{2\ve}})],
\end{align}
which by the independence of $\{\zeta_n\}$ implies that
\begin{align} \label{Cor-noniid-even}
&\sigma_{k_i,k_j}(1+h,1) \nonumber  \\
 =&4 a^{k_i+k_j-2} Var(\xi)
\sum_{\substack{(l,\ora{m}_l,\ora{n}_l) \in \Psi_{k_i} \\
(l',\ora{m}'_{l'},\ora{n}'_{l'}) \in \Psi_{k_j}}} C_l^{\ora{m}_l,
\ora{n}_l} C_{l'}^{\ora{m}'_{l'}, \ora{n}'_{l'}}
\sum_{\substack{1\leq q\leq l \\ 1\leq q' \leq l'}} m_q m'_{q'}
\delta_{q+h,q'}.
\end{align}

We shall compute explicitly the covariance $\Lambda(i,j)$ in
\eqref{Cov-nonstat.2}. For convenience, we place the leftmost vertex
of the circuit $\pi'$ of type $(l',\ora{m}'_{l'},\ora{n}'_{l'})$ at
$1$. Then, the circuit $\pi$ of the type
$(l,\ora{m}_{l},\ora{n}_{l})$ contributing to
\eqref{Cor-noniid-even} is the one with the leftmost vertex $1+h$,
and $\sum_{1\leq q\leq l} \sum_{1\leq q' \leq l'} m_q m'_{q'}
\delta_{q+h,q'} = 1/4 \sum_{q\geq 0}
\mathbb{I}_{h+q}(\pi)\mathbb{I}_{h+q}(\pi')$, where
$\mathbb{I}_i(\gamma)$ dentes the number of edges $\{i,i+1\}$ in the
circuit $\gamma$, i.e. $\mathbb{I}_i(\gamma)=2m_{i+1-j}$ if the
leftmost vertex of $\gamma$ is $j$ (see also \cite[p.187]{P09}).
Hence, the summation in $\sigma_{k_i,k_j}(1+h,1)$ consists of the
overlapped parts of all these two circuits $\pi$ and $\pi'$. Similar
argument holds for $\sigma_{k_i,k_j}(1,1+h)$, with the slight
modification that the leftmost vertex of $\pi'$ is $1-h$ (we here
extend the positions of vertices to the negative integers). Thus,
using the notation $\Gamma(k,l)$ in \cite[(3.1)]{P09}, we deduce
that the summation in \eqref{Cov-nonstat.2} is equal to $ \frac 14
\sum_{(\gamma_1,\gamma_2)\in \Gamma(k_i,k_j)} \sum_{i<0}
 \mathbb{I}_i(\gamma_1) \mathbb{I}_i(\gamma_2)
=\frac 14 k_ik_j
 {k_i \choose k_i/2} {k_j \choose  k_j/2}$
(see \cite[p.213]{P09} for the equality), which together with the
coefficients $(\a(k_i+k_j)+1-2\ve)^{-1}$ and
$4a^{k_i+k_j-2}Var(\xi)$ in \eqref{Cov-nonstat.2} and
\eqref{Cor-noniid-even} respectively yields
\eqref{Cor-noniid-covariance} when $k_i,k_j$ are even.

When $k$ is odd, \eqref{Qlmn.2} implies that $|\ora{n}_l|$ shall be
also odd. As in the preceding arguments,
\eqref{cal-abd.1-nonsym}-\eqref{cal-abd.3-nonsym} yield that the
only case contributing to the main order is that $|\ora{n}_l| =1$,
namely, there is only one loop in the circuit. We use the $q$-th
normal basis $\ora{e}_q$ of $\bbr^{[\frac k2]+1}$ to indicate that
the vertex of this loop is $n+q$ if the leftmost vertex of the
circuit is $n$, $0\leq q \leq [\frac k2]$. Then,
\eqref{cal-abd.1-nonsym} yields that
\begin{align} \label{Cor-X-odd}
   n^{-\a k } X_{k,n}
   = \sum_{\substack{(l,\ora{m}_l,\ora{e}_q)\in \Psi_k \\  0\leq q\leq [\frac{k}{2}]} }
    C_l^{\ora{m}_l,\ora{e}_q} a^{k-1} \zeta_{n+q} n^{-\ve}
    +\mathcal{O}(n^{-2\ve}).
\end{align}
(Note that, this includes the easy case $k=1$.) Hence, it follows
that
\begin{align} \label{Cor-noniid-odd}
 &\sigma_{k_i,k_j}(1+h,1) \nonumber  \\
 =& a^{k_i+k_j-2} Var(\zeta)
   \sum_{\substack{(l,\ora{m}_l,\ora{e}_q)\in \Psi_{k_i} \\ 0\leq q \leq  [\frac{k}{2}] }}
   \sum_{\substack{(l',\ora{m}'_l,\ora{e}_{q'})\in \Psi_{k_j} \\0\leq q' \leq [\frac{k}{2}]}}
    C_l^{\ora{m}_l,\ora{e}_q}
    C_{l'}^{\ora{m}'_{l'},\ora{e}_{q'}}
    \delta_{q+h,q'}.
\end{align}
In particular, $\sigma_{k_i,k_j}(1+h,1)$ vanishes when $0<\ve< \a$,
since in that case $\zeta =0$.

The computation of $\Lambda(i,j)$ in this case is easier. Similarly,
place the leftmost vertex of $\pi'$ of the type
$(l',\ora{m}'_l,\ora{n}'_l)$ at $1$. \eqref{Cor-noniid-odd}
indicates that leftmost vertex of $\pi$ of the type $(l,\ora{m}_l,
\ora{n}_l)$ is $1+h$, and the summation in \eqref{Cor-noniid-odd}
consists of all these two circuits $\pi$ and $\pi'$ with exactly one
loop at the same vertex. Thus, the summation in
\eqref{Cov-nonstat.2} is equal to $|\Gamma(k_i,k_j)|= k_ik_j {k_i-1
\choose  (k_i-1)/2} {k_j-1 \choose (k_j-1)/2}$ (see \cite[Remark
4]{P09} for this equality), which yields consequently
\eqref{Cor-noniid-covariance}.

Finally, as all $\xi_n$ and $\zeta_n$ are independent, it follows
from \eqref{Cor-X-even} and \eqref{Cor-X-odd} that
$\sigma_{k_i,k_j}(1+h,1)=0$ in the other cases, which implies
immediately
$\Lambda(i,j)=0$. The proof is complete. \hfill $\square$ \\

{\it $\beta$-Hermite ensemble (\cite{DE02,DE06}).} This model has
been widely studied in the literature. In this model, $a_i=b_i$,
$a_i$ is distributed as $\beta^{-1/2}\chi_{i\beta}$ ($\chi_{i\beta}$
is the $\chi$ distribution with $i\beta$ degrees of freedom),
$\beta>0$, $d_i$ is normally distributed as $N(0,2/\beta)$, and
$\{a_i,d_i\}$ are all independent. One significant fact is that the
density of the eigenvalue distribution is given by $C_{n\beta}
\prod_{i<j}
     |\lbb_i-\lbb_j|^{\beta}\ exp(-\beta  \sum_{i=1}^n
     \lbb_i^2 /2)$ with $C_{n\beta}$ a normalization. Hence, it
generalizes $\beta=1, 2, 4$ in  the classical Gaussian ensembles
(i.e. GOE, GUE, GSE respectively) to continuous exponents $\beta>0$,
which are connected to lattice gas theory, Selberg-type integrals,
Jack polynomials and so on. Now, using the fact that $\chi_r -
\sqrt{r} \overset{d}{\rightarrow} N(0,1/2)$, we deduce that in
Assumption $(H.2)$ $\a=\ve=1/2$, $a=1$, $\eta= N(0,1/(2\beta))$,
$d=0$ and $\zeta=N(0,2/\beta)$. Thus, Corollary \eqref{COR-NONIID}
implies Gaussian fluctuations of the traces and the covariances are
given by
\begin{align*}
   \Lambda(i,j) =
   \left\{
     \begin{array}{ll}
       \frac {1}{\beta} \frac{ k_ik_j}{k_i+k_j}
       {k_i \choose k_i/2}{k_j \choose k_j/2}, & \hbox{if $k_i, k_j$ even;} \\
       \frac {4}{\beta} \frac{k_ik_j}{k_i+k_j}
       {k_i -1 \choose (k_i -1)/2}{k_j -1 \choose (k_j -1)/2}, & \hbox{if $\ve=\a$, $k_i, k_j$ odd ;} \\
       0, & \hbox{otherwise,}
     \end{array}
   \right.
\end{align*}
which coincides with \cite[Corollary 2]{P09}\footnote{In
\cite[Corollary 2]{P09}, the terms $\frac{kl}{\a(k+l)}$ in the even
and odd cases shall be modified by
$\frac{kl\sigma_Z^2}{\a(k+l)+1-2\ve}$ and
$\frac{kl\sigma_d^2}{\a(k+l)+1-2\a}$ respectively.} and
\cite[Theorem 1.2]{DE06}.

\section{Deviations} \label{LDP-MDP}

This section is devoted to the large deviation and moderate
deviation principles for the traces. Due to technical reasons, we
will consider the tridiagonal random matrix \eqref{TRM} with entries
satisfying Assumption $(H.3)$ in the case $\a=0$ below.

\begin{enumerate}
\item[(H.3)] $(i)$. In the non-symmetric case, $\{(a_{i-1},d_i,b_i)\}$  are i.i.d. random
vectors. Set $\nu := \bbp \circ (a_1,d_2,b_2)^{-1} \in
\mathscr{P}(\bbr ^3)$.

In the symmetric case $(a_i=b_i, \forall i\geq 1)$, $\{a_i\}$ is a
sequence of i.i.d. random variables, and $\{d_i\}$ satisfies:

$(ii)$. $\{d_i\}$ is a sequence of i.i.d. random variables and
independent of $\{a_i\}$. Set $\nu':= \bbp \circ (d_1,a_1)^{-1} \in
\mathscr{P}(\bbr ^2)$. Or,

$(iii)$. $d_i=f(a_{i-1},a_i)$, where $f$ is a continuous function
and $a_0=0$. Set $\nu''=\bbp \circ (a_1)^{-1} \in
\mathscr{P}(\bbr)$.

Moreover, in all cases $\{a_i,d_i,b_i\}$ are bounded random
variables.
\end{enumerate}

\begin{remark}
Assumption $(H.3)$ already includes the models corresponding to the
case $\a=0$ in Subsection \ref{CLT-APP}, such as the Anderson model,
the Hatano-Nelson model, the random birth-death Markov kernel and
the random birth-death $Q$ matrix.
\end{remark}

\subsection{Large deviations} \label{LDP}

Let us first introduce some basic notations and definitions for
large deviation principles. More details can be found in \cite{DZ98}
and \cite{E88}. Then, we formulate and prove our result in Theorem
\ref{THM-LDP} below.

Consider a complete separable metric space $\mathscr{X}$,
$\mathscr{P}(\mathscr{X})$ denotes all the probability measures in
$\mathscr{X}$. A sequence of probability measures $\{\mu_n\}\subset
\mathscr{P}(\mathscr{X})$ is said to satisfy the large deviation
principle (LDP for short) with speed $s_n \to \9$ and good rate
function $I: \mathscr{X} \to [0,\9]$, if the level sets $\{x \in
\mathscr{X}: I(x)\leq c\}$ are compact for all $c\in[0,\9)$ and if
for all Borel set $A$ of $\mathscr{X}$,
\begin{align*}
  -\inf\limits_{x\in A^o} I(x)
  \leq& \liminf\limits_{n\to\9} \frac{1}{s_n} \log \mu_n (
  A)
  \leq \limsup\limits_{n\to\9} \frac{1}{s_n} \log \mu_n (
  A)
  \leq - \inf\limits_{x\in \overline{A}} I(x),
\end{align*}
where $A^o$ and $ \overline{A}$ denote the interior and closure of
$A$ respectively. In that case, we shall simply say that $\{\mu_n\}$
satisfies the $LDP(s_n, I)$ on $\mathscr{X}$. We also say that a
family of $\mathscr{X}$-valued random variables satisfies the
$LDP(s_n, I)$ if the family of their laws does. For any $\mu,\mu'
\in \mathscr{P}(\mathscr{X})$,  the relative entropy of $\mu$ with
respect to $\mu'$ is defined by
\begin{align*}
   H (\mu |\mu') =\left\{
                    \begin{array}{ll}
                      \int g \log g d\mu', & \hbox{if $g= \frac{d \mu}{d\mu'}$ exists;} \\
                      \9, & \hbox{otherwise.}
                    \end{array}
                  \right.
\end{align*}

We also consider the product space $\mathscr{X}^r$ with $r\geq 2$. A
probability measure $\mu\in \mathscr{P}(\mathscr{X}^r)$ is said to
be shift invariant if for any Borel set $A$ of $\mathscr{X}^{r-1}$,
\begin{align*}
   \mu (x\in \mathscr{X}^r: (x_1,\cdots,x_{r-1}) \in A)
   = \mu (x\in \mathscr{X}^r: (x_2,\cdots,x_{r}) \in A).
\end{align*}
Moreover, for any $\mu\in \mathscr{P}(\mathscr{X}^{r-1})$, $\rho \in
\mathscr{P}(\mathscr{X})$, define the probability measure $\mu
\otimes_{r} \rho \in \scrp(\scrx^r)$ by, for any Borel set $A$ of
$\scrx^r$,
\begin{align*}
   (\mu \otimes_{r} \rho)\ (A)
   =\int\limits_{\scrx^{r-1}} \mu(dx) \int\limits_{\scrx} I_{\{(x,y)\in
   A\}} \rho(dy),
\end{align*}
where $I_{\{(x,y)\in
   A\}}$ is the indicator function of the set $\{(x,y)\in
   A\}$.  Define for $\rho \in \scrp(\scrx)$,
\begin{align} \label{rate-LDP}
   I_{r,\rho}(\mu)
   =\left\{
      \begin{array}{ll}
        H(\mu | \mu_{r-1} \otimes_{r} \rho), & \hbox{if $\mu$ is shift invariant;} \\
        \9, & \hbox{otherwise,}
      \end{array}
    \right.
\end{align}
where $\mu_{r-1} $ denotes the marginal of $\mu$ on the first
$(r-1)$ coordinates.

Now, come to the tridiagonal random matrix \eqref{TRM}. For each
$k\geq 1$, set $r_k=[\frac k2]+1$. In the non-symmetric case, take
$\scrx= \bbr^3$, and define the function $F$ on $(\bbr^3)^{r_k}$ by
\begin{align} \label{Xi-nonsym}
  F(\a_1,\cdots, \a_{r_k})
  = \sum\limits_{(l,\ora{m}_l,\ora{n}_l)\in \Psi_k}
    C_l^{\ora{m}_l,\ora{n}_l}
    \prod\limits_{j=0}^{l-1} (x_{1+j}z_{1+j})^{m_{j+1}}
    \prod\limits_{j=0}^l y_{1+j}^{n_j},
\end{align}
where $\a_i = (x_{i-1},y_i,z_i)\in \bbr^3$, $1\leq i \leq r_k$ and
$x_0=0$. Associating with $F$, define $I_{3r_k, \nu}^{F}$ on
$\bbr^3$ by
\begin{align*}
   I_{3r_k,\nu}^{F} (x) := \inf \{ I_{r_k,\nu}(\mu):\ \mu\in\mathscr{P}((\bbr^3)^{r_k}),\  x=\mu
   (F)\},
\end{align*}
where $\nu$ is defined as in Assumption $(H.3)$, and $\mu(F)$
denotes the integration of $F$ with respect to $\mu$.

Similarly, in the symmetric case $(ii)$ in $(H.3)$, take $\scrx =
\bbr^2$, define the function $F'$ on $(\bbr^2)^{r_k}$ by
\begin{align} \label{Xi-sym}
  F'(\beta_1,\cdots, \beta_{r_k})
  = \sum\limits_{(l,\ora{m}_l,\ora{n}_l)\in \Psi_k}
    C_l^{\ora{m}_l,\ora{n}_l}
    \prod\limits_{j=0}^{l-1} z_{1+j}^{2m_{j+1}}
    \prod\limits_{j=0}^l y_{1+j}^{n_j}
\end{align}
with $\beta_{i}=(y_i,z_i)$, $1\leq i \leq r_k$, and the
corresponding function $I_{2r_k, \nu'}^{F'}$ is defined by
\begin{align*}
   I_{2r_k,\nu'}^{F'} (x):= \inf \{ I_{r_k,\nu'}(\mu):\ \mu\in\mathscr{P}((\bbr^{2})^{r_k}),\ x= \mu (F')
\}
\end{align*}
where $\nu'$ is defined as in Assumption $(H.3)$.

In the symmetric case $(iii)$ in $(H.3)$, take $\scrx=\bbr$, and
define the function $F''$ on $\bbr^{r_k+1}$ by
\begin{align*}
    F''(z_0,\cdots,z_{r_k}) = F' (\gamma_1, \cdots, \gamma_{r_k}),
\end{align*}
where $\gamma_i= (f(z_{i-1}, z_i), z_i)$, $1\leq i \leq r_k$. The
related function $I_{r_k+1, \nu''}^{F''}$ is defined  analogously by
\begin{align*}
   I_{r_k+1,\nu''}^{F''} (x):= \inf \{ I_{r_k+1,\nu''}(\mu):\
\mu\in\mathscr{P}((\bbr)^{r_k+1}),\ x= \mu (F'') \}
\end{align*}
with $\nu''$ defined as in $(iii)$ in $(H.3)$.

We are ready to state our large deviation result. Since when $k=1$,
$Tr Q_n^k = \sum_{i=1}^n d_i$, i.e. the sum of i.i.d random
variables, of which the large deviation is well known by the
Cram\'{e}r theorem (see \cite[Theorem 2.2.3]{DZ98}). Therefore,
below we are concerned with the case $k\geq 2$.

\begin{theorem} \label{THM-LDP}
Assume $(H.3)$. Let $k\geq 2$ and set\ $r_k:= [\frac k2]+1$. Then,
$\{\frac 1n Tr Q_n^k\}_{n\geq 1}$ satisfies the
$LDP(n,I_{3r_k,\nu}^{F})$, $LDP(n,I_{2r_k,\nu'}^{F'})$ and
$LDP(n,I_{r_k+1,\nu''}^{F''})$ in the cases $(i)-(iii)$ in $(H.3)$
respectively.
\end{theorem}

{\it Proof. } Let us first consider the non-symmetric case $(i)$ in
$(H.3)$. We note that $\{\frac 1n Tr Q_n^k\}$ and $\{\frac 1n
\sum_{i=1}^n X_{k,i} \}$ is exponentially equivalent, i.e., for any
$\delta
>0$,
\begin{align} \label{equi-LDP-TrQk}
   \lim\limits_{n\to\9} \frac 1n \log \bbp(\frac 1n |Tr Q_n^k - \sum\limits_{i=1}^n X_{k,i}|> \delta)
   =-\9.
\end{align}
(See the Appendix for the proof.) Thus, by the exponential
equivalence theorem(\cite[Theorem 4.2.13]{DZ98}), the proof reduces
to proving the large deviation of $ \frac 1n \sum_{i=1}^n X_{k,i}$.

Now, set the random vectors $\a_i=(a_{i-1},d_i,b_i) \in \bbr^3$,
$i\geq 1$. By the independence and the identical distribution of
$\{\a_i\}$ in Assumption $(H.3)$, we regard $\{\a_i\}$ as a Markov
chain on $\bbr^3$ with the transition probability $\pi(x,dy) : =
\nu(dy)$, $x,y\in \bbr^3$. Note that $\{\pi(x,dy)\}$ obviously
satisfy the uniform Assumption $(U)$ in \cite[p.275]{DZ98} (see also
Hypothesis $1.1 (a)$ in \cite{E88}), i.e., $\{\a_i\}$ is a uniform
Markov chain on $\bbr^3$. Thus, applying \cite[Theorem 6.5.12]{DZ98}
(see also \cite[Theorem 1.4]{E88}), we have that the multivariate
empirical measure $\mu_n:=\frac 1n \sum_{i=1}^n
\delta_{(\a_{i},\cdots,\a_{i+r_k-1})}$ satisfies the $LDP(n,I_{
r_k,\nu})$ with the good rate function $I_{r_k, \nu}$ defined as in
\eqref{rate-LDP}.

On the other hand, by the definition of $X_{k,i}$ in \eqref{Xi}, we
note that
\begin{align} \label{X-F-nonsym}
X_{k,i} = F(\a_{i}, \cdots, \a_{i+r_k-1}),
\end{align} namely, $X_{k,i}$ can be viewed as a continuous
function on the product space $(\bbr^3)^{r_k}$. Moreover,
\begin{align*}
\frac 1n \sum\limits_{i=1}^n X_{k,i} = \frac 1n \sum\limits_{i=1}^n
F(\a_{i},\cdots, \a_{i+ r_k-1})= \mu_n(F),
\end{align*}
which implies that $\frac 1n \sum_{i=1}^n X_{k,i}$ is an additive
functional of the uniform Markov chain.

Therefore, we can apply the contraction principle (see \cite[Theorem
4.2.1]{DZ98}) to obtain the large deviation of $\frac 1n
\sum_{i=1}^n X_{k,i}$ as specified in the non-symmetric case.

The symmetric case can be treated analogously. In fact, in the case
$(ii)$ in $(H.3)$, $X_{k,i} = F'(\beta_i,\cdots, \beta_{i+r_k-1})$,
and the random vectors $\beta_i = (d_i,b_i)\in \bbr^2$ forms a
uniform Markov chain on $\bbr^2$ with the transition probability
$\pi'(x,dy)=\nu'(dy)$. Moreover, in the case $(iii)$ in $(H.3)$,
$X_{k,i} = F''(b_{i-1},\cdots, b_{i+r_k-1})$, and $\{b_i\}\subset
\bbr$ forms a uniform Markov chain in $\bbr$ with
$\pi''(x,dy)=\nu''(dy)$. Therefore, applying \cite[Theorem
6.5.12]{DZ98} and the contraction principle, we obtain the asserted
large deviation results. The proof of Theorem \ref{THM-LDP} is
complete. \hfill $\square$

\subsection{Moderate deviations} \label{MDP}

Moderate deviation principles for dependent random variables are
widely studied in the literature, see e.g. \cite{W95, DMPU09, MPR11}
and references therein. Here, for the moderate deviations of the
traces, we prefer to give an elementary proof based on the blocking
arguments as those in the proof of Theorem \ref{THM-CLT-NONSTA}.

\begin{theorem} \label{THM-MDP-H1}
Assume $(H.3)$. Set $\lbb_n=n^{-\nu}$ with $\nu\in (0,1)$, and let
$Tr \wt{Q_n^k} = Tr Q_n^k - \bbe Tr Q_n^k$. Then, for every $k \geq
1$ and  $D_k>0$, where $D_k$ is defined as in \eqref{Cor.1-Dk},
$\{\sqrt{\frac {\lbb_n}{n}} Tr \wt{Q_n^k} \}_{n\geq 1}$ satisfies
the $LDP(\lbb_n^{-1}, \frac{x^2}{2D_k})$.
\end{theorem}

{\it Proof of Theorem \ref{THM-MDP-H1}.} Let $ \wt{X}_{k,i} =
X_{k,i} - \bbe X_{k,i}$. We first note that $\{ \sqrt{\frac
{\lbb_n}{n}} Tr \wt{Q_n^k} \}$ and $\{\sqrt{\frac {\lbb_n}{n}}
\sum_{i=1}^n \wt{X}_{k,i}\}$ is exponentially equivalent, i.e. for
any $\delta >0$,
\begin{align} \label{equi-MDP-TrQk}
   \limsup\limits_{n\to \9} \lbb_n \log \mathbb{P}
   \(\sqrt{\frac{\lbb_n}{n}} |Tr \wt{Q_n^k} - \sum\limits_{i=1}^n \wt{X}_{k,i}| \geq
    \delta \) = -\9.
\end{align}
(The proof is postponed to the Appendix.) Hence, it is equivalent to
prove the moderate deviation for $\sum_{i=1}^n \wt{X}_{k,i}$. For
this purpose, we will use the blocking arguments as in the proof of
Theorem \ref{THM-CLT-NONSTA}.

Let $p=n^{\nu + \ve}$, where $\ve \in (\frac 34 (1-\nu),1-\nu)$. Set
$n'=[\frac{n}{p}]$ and $r=n-pn'$. Let $\wt{Y}_{n,i}, \wt{U}_{n,i},
\wt{Z}_{n,i}$ and $\wt{T}_n$ be as in the proof of Theorem
\ref{THM-CLT-NONSTA}, but with $\a=0$. Then
\begin{align*}
   \sqrt{\frac{\lbb_n}{n}} \sum\limits_{i=1}^n \wt{X}_{k,i}
   =  \sqrt{\frac{\lbb_n}{n}}\sum\limits_{i=1}^p \wt{U}_{n,i}
      + \sqrt{\frac{\lbb_n}{n}}   \wt{T}_n.
\end{align*}
Moreover, denote by $\Lambda_n$ the logarithmic moment generating
function of $\sqrt{\frac {\lbb_n}{n}} \sum\limits_{i=1}^p
\wt{U}_{n,i}$, i.e. $\Lambda_n(t) = \log \bbe \exp( t \sqrt{\frac
{\lbb_n}{n}} \sum_{i=1}^p \wt{U}_{n,i})$, $t\in \mathbb{R}$.

We shall prove below that, for every $\delta>0$,
\begin{align} \label{equi-MDP-Ui}
    \limsup \limits_{n\to \9} \lbb_n \log \bbp \(|\sqrt{\frac{\lbb_n}{n}} \wt{T}_n| \geq \delta \) = -\9,
\end{align}
and
\begin{align} \label{MDP-Ui}
    \lim\limits_{n\to \9} \lbb_n  \Lambda_n (\lbb_n^{-1} t) = \frac {t^2}{2} D_k,
\end{align}
where $D_k$ is the variance defined as in \eqref{Cor.1-Dk}. Then, by
the exponential equivalence, \eqref{equi-MDP-Ui} implies that we
only need to consider the moderate deviations of $ \sum_{i=1}^p
\wt{U}_{n,i}$. Consequently, \eqref{MDP-Ui} and the
G\"{a}rtner-Ellis theorem (see e.g. \cite{DZ98})
yield the asserted moderate deviation principle for the traces.\\

It remains to prove \eqref{equi-MDP-Ui} and  \eqref{MDP-Ui}. For the
proof of \eqref{equi-MDP-Ui}, since $\wt{T}_n = \sum_{i=1}^p
\wt{Z}_{n,i}$ and $\wt{Z}_{n,i}$ are independent, using Assumption
$(H.3)$ and the Bernstein inequality (see e.g. \cite[p.21]{BS10}),
we have
\begin{align*}
    \bbp (|\sqrt{\frac{\lbb_n}{n}} \wt{T}_n| \geq \delta)
    \leq 2 e^{- \delta^2/(2(B_n^2+ c\delta))},
\end{align*}
where $ B_n^2 = \frac{\lbb_n}{n} \bbe |\wt{T}_n|^2 =
\mathcal{O}(\frac{\lbb_n p }{n})$ and
$c=\sqrt{\frac{\lbb_n}{n}}\sup\limits_{1\leq i\leq
p}\|\wt{Z}_{n,i}\|_{\9}=\mathcal{O}(\sqrt{\frac{\lbb_n}{n}})$.
Thus,
\begin{align*}
     \lbb_n \log \bbp (|\sqrt{\frac{\lbb_n}{n}} \wt{T}_n| \geq \delta)
     \leq& \lbb_n \log (2 e^{- \delta^2/(2(B_n^2+ c\delta))}) \\
     =& \lbb_n \log 2 - \frac{\delta^2}{\mathcal{O}(\frac{p}{n} + \frac{1}{\sqrt{\lbb_n
     n}})}  \to -\9,
\end{align*}
which implies \eqref{equi-MDP-Ui}, as claimed. \\

Coming to the proof of \eqref{MDP-Ui}, we note that by the uniform
boundedness in Assumption $(H.3)$ and $\ve > \frac 12 (1-\nu)$,
$$\frac{1}{\sqrt{\lbb_n n}}
     \wt{U}_{n,i} = \mathcal{O}(n'/\sqrt{\lbb_n n})= \mathcal{O}(n^{\frac 12(1 -
     \nu)
     -\ve})=o(1).$$
Hence, it follows that
\begin{align} \label{MDP-Log.1}
    \lbb_n \Lambda_n(\lbb_n^{-1} t)
    =&  \sum\limits_{i=1}^p \lbb_n \log \bbe \exp( \frac{t}{\sqrt{\lbb_n n}}
     \wt{U}_{n,i}) \nonumber \\
    =& \sum\limits_{i=1}^p \lbb_n \log \left[1+ \frac{t^2}{2 \lbb_n n } \bbe  \wt{U}_{n,i}^2 + \mathcal{O}(\frac{\bbe
     \wt{U}_{n,i}^3}{(\lbb_n n)^{\frac{3}{2}}}) \right].
\end{align}
Similarly to \eqref{exp-uni} and \eqref{cal-u}, we have
\begin{align*}
   \frac{1}{\lbb_n n} \bbe \wt{U}_{n,i}^2
   =&  \frac{1}{\lbb_n n} \left[\bbe \(\sum\limits_{h=1}^{m_k} \wt{Y}_{n,(i-1)n'+h} \)^2
      + \sum\limits_{h=1}^{n'-2m_k} A_{n,(i-1)n'+h} \right]\\
   =& \mathcal{O}(\frac{1}{\lbb_n n}) +  \frac{1}{\lbb_n n} \sum\limits_{h=1}^{n'-2m_k}
   A_{n,(i-1)n'+h},
\end{align*}
where $A_{n,(i-1)n'+h}$ are defined as in \eqref{def-ani}. Hence, it
follows from Assumption $(H.3)$ and similar computations as in
\eqref{cal-u}-\eqref{esti-2nd} that
\begin{align} \label{MDP-esti-2nd}
   \frac{1}{\lbb_n n} \bbe \wt{U}_{n,i}^2
   = \mathcal{O}(\frac{1}{\lbb_n n} )+ \frac{1}{\lbb_n p }(D_k +o(1)).
\end{align}
Moreover, for the third moment in \eqref{MDP-Log.1}, by the choices
of $\lbb_n$ and $n'$,
\begin{align} \label{MDP-esti-3rd}
   \frac{1}{(\lbb_n n)^{\frac{3}{2}}} \bbe \wt{U}_{n,i}^{3}
   =& \frac{1}{(\lbb_n n)^{\frac{3}{2}}} \bbe (\sum\limits_{j=(i-1)n'+1}^{in'-m_k}
   \wt{Y}_{n,j})^3  \nonumber \\
   =& \frac{1}{(\lbb_n n)^{\frac{3}{2}}} \mathcal{O}((n')^3) \nonumber  \\
   =& n^{3(\frac{1}{2} - \frac 12 \nu - \ve)}
   =o(1).
\end{align}
Consequently, plugging  \eqref{MDP-esti-2nd} and
\eqref{MDP-esti-3rd} into \eqref{MDP-Log.1}, we obtain that
\begin{align*}
    \lbb_n \Lambda_n (\lbb_n^{-1} t)
   =&\sum\limits_{i=1}^p  \lbb_n \left[\frac{t^2}{2\lbb_n n} \bbe \wt{U}_{n,i}^2
     + \mathcal{O}(\frac{\bbe \wt{U}_{n,i}^3}{(\lbb_n n)^{\frac 32}})\right]  +o(1)  \\
   =& \mathcal{O}(\frac{1}{n'}) + \frac{t^2}{2} D_k + \mathcal{O}(n^{\frac{3}{2} - \frac 32 \nu -
   2\ve}) + o(1) \\
   \to&  \frac{t^2}{2}D_k.
\end{align*}
which implies \eqref{MDP-Ui}, thereby completing the proof. \hfill
$\square$\\

\section{Discussions} \label{DIS}

$1$. For more general finite diagonal random matrix (see e.g.
\cite{BS13}, \cite{P09}), the Gaussian fluctuations and deviations
of the traces can be treated in a similar way. Indeed, using the
finite width of band, we can still reduce the asymptotical analysis
of the traces to those of the corresponding $m$-dependent random
variables.

$2$. It seems difficult to prove deviation results for the general
case $\a>0$. In some special cases, e.g. the GUE case ($\a=1/2$),
deviation results are known for the empirical spectral distribution
(hence for the traces by the contraction principle), we refer to
\cite{AGZ10, DGZ03}.

$3$. In the derivation of the large deviation results in Subsection
\ref{LDP}, we regard the trace as an additive functional of a
uniform Markov chain. With this point of view, one can expect to
achieve similar results for the traces when the entries of
\eqref{TRM} form an appropriate Markov chain, e.g. positive Harris
recurrent. We refer the interested reader to \cite{AN14} and
references therein.

\section{Appendix} \label{APP}

{\it Proof of \eqref{TrQ-SumX}.} By Assumption $(H.1)$ and $(H.2)$,
\begin{align} \label{bdd-Q}
  \sup\limits_{i\geq 1} \bbe |n^{-\a k} \wt{Q}_{l,i}^{\ora{m}_l,
  \ora{n}_l}|^2 < \9,
\end{align}
where $\wt{Q}_{l,i}^{\ora{m}_l} = Q_{l,i}^{\ora{m}_l} - \bbe
Q_{l,i}^{\ora{m}_l}$. Hence,
\begin{align*}
   &\|n^{-\a k} (Tr \wt{Q_n^k} - \sum\limits_{i=1}^n \wt{X}_{k,i}) \|_2 \\
   =&\|\sum\limits_{(l,\ora{m}_l,\ora{n}_l)\in \Psi_k} C_l^{\ora{m}_l, \ora{n}_l}
   \sum\limits_{i=n-l+1}^n n^{-\a k} \wt{Q}_{l,i}^{\ora{m}_l,
   \ora{n}_l}\|_2 \\
   \leq& \sum\limits_{(l,\ora{m}_l,\ora{n}_l)\in \Psi_k} C_l^{\ora{m}_l, \ora{n}_l}
   \sum\limits_{i=n-l+1}^n \| n^{-\a k} \wt{Q}_{l,i}^{\ora{m}_l,
   \ora{n}_l}\|_2
   =\mathcal{O}(1),
\end{align*}
with $\|\cdot\|_2$ denoting the standard $L^2$ norm, which implies
\eqref{TrQ-SumX}. \hfill $\square$\\

{\it Proof of \eqref{asym-cov}.}  It is equivalent to prove that, as
$n\to \9$,
\begin{align} \label{asym-XXj}
   n^{-2\a k} Cov (X_{k,n}, X_{k,n+j})
    \to Cov (Z_{k,1}, Z_{k,1+j}).
\end{align}
First consider the nonsymmetic case. Let $\a_i=(a_{i-1}, d_i, b_i)$,
$i\geq 1$ and $a_0=0$. By the independence and weak convergence of
$\a_i$ in Assumption $(H.1)$, it is not difficult to deduce that
\begin{align*}
   n^{-\a k} (\a_n, \cdots, \a_{n+[\frac k2]})
   \overset{d}{\rightarrow} (\wt{\a}_1, \cdots, \wt{\a}_{1+[\frac
   k2]}),
\end{align*}
where $\wt{\a}_i$ are independent but with the same distribution as
that of $(a+ \eta, d+\zeta, b+ \xi)$. Then, by \eqref{X-F-nonsym}
and the continuous mapping theorem (\cite[Theorem 3.2.4]{D05}),
\begin{align} \label{converg-X}
   n^{-\a k} X_{k,n} \overset{d}{\rightarrow} F(\wt{\a}_1, \cdots, \wt{\a}_{1+[\frac{k}{2}]})
\end{align}
with $F$ the continuous function defined as in \eqref{Xi-nonsym}.

Similarly,
\begin{align} \label{converg-Xj}
   n^{-\a k} X_{k,n+j} \overset{d}{\rightarrow} F(\wt{\a}_{1+j}, \cdots,
   \wt{\a}_{1+j+[\frac{k}{2}]}).
\end{align}

On the other hand, by \eqref{asym-eatazeta.2} and H\"{o}lder's
inequality
\begin{align} \label{uninteg-X4}
   \sup\limits_{n\geq 1} \bbe (n^{-\a k} X_{k,n})^4 < \9,
\end{align}
which implies the uniform integrabilities of $n^{-2\a k}
X_{k,n}X_{k,n+j}$, $n^{-\a k} X_{k,n}$ and $n^{-\a k} X_{k,n+j}$,
$n\geq 1$.

Therefore, by \eqref{converg-X}-\eqref{uninteg-X4}, we can apply the
Skorohod representation theorem and the uniform integrability to
take the limit and obtain \eqref{asym-XXj} for the non-symmetric
case. (See e.g. \cite[Theorem 3.2.4]{D05} for similar arguments for
the bounded continuous mapping.)\\

The symmetric case can be proved analogously. In fact, with the
uniform integrability \eqref{uninteg-X4}, we only need to check the
weak convergence of $n^{-\a k} X_{k,n}$.

In the case that all $d_i$ are independent of $a_i(=b_i)$, let
$\beta_i = (d_i, a_i)$, $i\geq 1$. In this case, $X_{k,n}=
F'(\beta_n,\cdots, \beta_{n+[\frac k2]})$ with $F'$ defines as in
\eqref{Xi-sym}. Then, it follows from similar arguments as above
that
\begin{align}
   n^{-\a k} X_{k,n} \overset{d}{\rightarrow} F'(\wt{\beta_1}, \cdots, \wt{\beta}_{1+[\frac
   k2]}),
\end{align}
where $\wt{\beta_i}$, $1\leq i\leq [\frac k2]+1$, are independent
but with the common distribution as that of $(d+\xi, a+\eta)$.

In the case that $d_i=f(a_{i-1}, a_i)$, then $\beta_i = (f(a_{i-1},
a_i), a_i)$, $X_{k,n} $ is now a continuous function of $a_{n-1},
\cdots, a_{n+[\frac k2]}$. As $(a_{n-1},\cdots, a_{n+[\frac k2]})
\overset{d}{\rightarrow} (\wt{\a_1}, \cdots, \wt{\a}_{[\frac
k2]+2})$, by the continuous mapping theorem, we can obtain the weak
convergence of $n^{-\a k} X_{k,n}$. The proof is consequently
complete. \hfill
$\square$ \\

{\it Proof of \eqref{equi-LDP-TrQk}.} Set $N(k)= |\Psi_k|$, the
number of sets in $\Psi_k$ which is defined as in \eqref{Psi}.
$N(k)$ is finite and depends only on $k$. Note that
\begin{align*}
   & \bbp \( \frac 1n  \left|Tr Q_n^k - \sum\limits_{i=1}^n X_i \right| \geq
   \delta \)  \\
   \leq & \sum\limits_{(l,\ora{m}_l,\ora{n}_l)\in \Psi_k}
        \bbp \( \left| C_{l}^{\ora{m}_l,\ora{n}_l} \sum\limits_{i=n-l+1}^n Q_{l,i}^{\ora{m}_l,\ora{n}_l} \right| \geq
    \frac{n \delta}{ N(k)} \).
\end{align*}
Then, letting $C$ denote the maximum of
$C_{l}^{\ora{m}_l,\ora{n}_l}$ over the finite sets in $\Psi_k$ and
setting $\delta_{k} = \delta/([\frac k2] C N(k))$, we deduce that
\begin{align*}
   \bbp \(\frac 1n \left|Tr Q_n^k - \sum\limits_{i=1}^n X_{k,i} \right| \geq
   \delta \)
   \leq \sum\limits_{(l,\ora{m}_l,\ora{n}_l)\in \Psi_k}
          \sum\limits_{i=n-l+1}^n
    \bbp \( \left|Q_{l,i}^{\ora{m}_l,\ora{n}_l} \right|
    \geq n \delta_{k} \).
\end{align*}
Since by \eqref{Qlmn} and Assumption $(H.3)$, $\{ |Q
_{l,i}^{\ora{m}_l,\ora{n}_l}| \}_{i\geq 1}$ is uniformly bounded,
which implies that  $\bbp (|Q _{l,i}^{\ora{m}_l,\ora{n}_l}| \geq
n\delta_k) =0$  for $n$ large enough.  Hence, by Lemma $1.2.15$ in
\cite{DZ98},
\begin{align} \label{esti-equi-LDP-TrQk}
    & \frac 1 n \log \mathbb{P} \(\frac 1n |Tr Q_n^k - \sum\limits_{i=1}^n X_{k,i}| \geq
   \delta \) \nonumber \\
  \leq& \max_{\substack{(l,\ora{m}_l,\ora{n}_l)\in \Psi_k \\ 1\leq j\leq [\frac k2]}}
     \frac 1 n \log \bbp \(|Q_{l, n-l+j}^{\ora{m}_l,\ora{n}_l}|
   \geq  n  \delta_{k} \)
  =-\9.
\end{align}
yielding \eqref{equi-LDP-TrQk} as claimed. \hfill $\square$\\

{\it Proof of \eqref{equi-MDP-TrQk}.} Similarly to the proof of
\eqref{esti-equi-LDP-TrQk}, we derive that
\begin{align*}
   &\lbb_n \log \mathbb{P} \(\sqrt{\frac{\lbb_n}{n}} |Tr \wt{Q_n^k} - \sum\limits_{i=1}^n \wt{X}_{k,i}| \geq
   \delta \) \\
  \leq& \max_{\substack{(l,\ora{m}_l,\ora{n}_l)\in \Psi_k \\ 1\leq j\leq [\frac k2]}}
     \lbb_n \log \bbp \(|\wt{Q}_{l, n-l+j}^{\ora{m}_l,\ora{n}_l}|
   \geq \sqrt{\frac{n}{\lbb_n}}  \delta_{k} \)
\end{align*}
with $\delta_k$ defined as in the proof of
\eqref{esti-equi-LDP-TrQk}, which yields \eqref{equi-MDP-TrQk}, due
to the fact that  $\{Q_{l, i}^{\ora{m}_l,\ora{n}_l}\}_{i\geq 1}$ is
uniformly bounded and
$n/\lbb_n \to \9$.  \\

{\it \bf Acknowledgement.} The author is grateful to Professor Dong
Han for several valuable discussions.

\end{document}